\documentclass[12pt,twoside]{amsart}

\usepackage{amssymb}

\newtheorem{thm}{Theorem}

\newtheorem{lem}[thm]{Lemma}
\newtheorem{cor}[thm]{Corollary}

\newtheorem{prop}[thm]{Proposition}

\theoremstyle{definition}
\newtheorem{defn}[thm]{Definition}

\newtheorem{say}[thm]{}

\newtheorem{rem}[thm]{Remark}          
\newtheorem{ack}{Acknowledgments}

\newtheorem{defn-thm}[thm]{Definition--Theorem}  

\theoremstyle{remark}
  
\newtheorem{case}{Case}

\renewcommand{\c}[0]{{\mathbb C}}  

\renewcommand{\o}[0]{{\mathcal O}} 
\newcommand{\z}[0]{{\mathbb Z}}

\renewcommand{\a}[0]{{\mathbb A}}

\newcommand{\p}[0]{{\mathbb P}}

\newcommand{\q}[0]{{\mathbb Q}}

\newcommand{\qtq}[1]{\quad\mbox{#1}\quad}

\newcommand{\diag}[0]{\operatorname{diag}}

\begin{document}
\bibliographystyle{amsplain}

\title[Fano hypersurfaces  in 
weighted projective 4-spaces]{Fano hypersurfaces  in \\
weighted projective 4-spaces}
\author{Jennifer M.\ Johnson   and J\'anos Koll\'ar}

\maketitle

A {\it Fano} variety is a projective  variety
whose  anticanonical class is ample.
A 2--dimensional Fano variety is called a Del Pezzo surface. 
 In higher dimensions, attention originally centered on smooth
 Fano 3--folds, but singular Fano varieties are also of
considerable interest in connection with the minimal model program.
The existence  of K\"ahler--Einstein metrics 
on Fano varieties has also been explored,
see \cite{bourg} for a summary of the main results. Here
 again the smooth case is of primary interest, but 
Fano varieties with quotient singularities and their 
orbifold metrics have also been studied.

In any given dimension there are only
finitely many families of smooth Fano varieties \cite{camp, nadel, kmm3},
but very little is known about them in dimensions 4 and up.
By allowing singularities, infinitely many families appear
and their distribution is very poorly understood.

A natural experimental testing ground is given by
hypersurfaces and complete intersections in 
weighted
projective spaces. These can be written down rather explicitly,
but they still provide many more examples than ordinary 
projective spaces. Experimental lists of certain
3--dimensional complete intersections were compiled by
\cite{fletch}. In connection with K\"ahler--Einstein metrics,
the 2--dimensional cases were first investigated
in  \cite{dk} and later in
\cite{jk1}. 

It is also of interest to study Calabi--Yau hypersurfaces 
and hypersurfaces  of general type in 
weighted
projective spaces. \cite{fletch} contains
some lists with terminal singularities.

The aim of this paper is threefold. 

First, we determine the complete
list of  anticanonically embedded quasi smooth 
Fano hypersurfaces in weighted projective 4-spaces.
There are 48 infinite series and 4442 sporadic examples (\ref{list.thm}).
As a consequence we obtain that the Reid--Fletcher list 
(cf.\ \cite[II.6.6]{fletch})
of 95  types  of anticanonically embedded quasi smooth 
terminal Fano threefolds
in  weighted projective 4-spaces is complete (\ref{95.cor}).

Second, we  prove that many of these Fano hypersurfaces 
admit a K\"ahler--Einstein metric (\ref{KE.cor}).
We also study the nonexistence of 
 tigers on these Fano 3--folds 
(in the colorful terminology of \cite{keel-mc}). 

Third, we prove that there are only finitely many families of
quasi smooth  
Calabi--Yau hypersurfaces in weighted projective spaces
of any given dimension (\ref{CY.thm}).
 This implies finiteness for various families of
general type hypersurfaces (\ref{gen.t.cor}).

\begin{defn}\label{wps.defn} For positive integers $a_i$
let $\p(a_0,\dots,a_n)$ denote the 
{\it  weighted
projective $n$-space} with weights $a_0,\dots,a_n$. 
(See \cite{dolg} or \cite{fletch}  for the basic definitions and results
on weighted
projective spaces.)
We always assume
that  any $n$  of  the
  $a_i$ are relatively prime. 
We frequently write $\p$ 
to denote a weighted
projective $n$-space
if the weights are irrelevant or clear from the context.
We use $x_0,\dots,x_n$ to denote the corresponding weighted
projective coordinates. 
$P_i\in \p(a_0,\dots,a_n)$ denotes the point
$(x_j=0\ \forall j\neq i)$.  
These are sometimes called the {\it vertices} of the weighted
projective space. (The vertices are uniquely determined if none of the
$a_i$  divides another.)
The affine chart where
$x_i\neq 0$  can be written as
$$
\c^n(y_0,\dots, \hat{y_i},\dots,y_n)/\z_{a_i}(a_0,\dots, \hat{a_i},\dots,a_n).
\eqno{(\ref{wps.defn}.1)}
$$
(Here and later  $\ \hat{}\ $ denotes an omitted coordinate.)
This shorthand denotes the quotient of $\c^n$ by the action
$$
(y_0,\dots, \hat{y_i},\dots,y_n)\mapsto 
(\epsilon^{a_0}y_0,\dots, \hat{y_i},\dots, \epsilon^{a_n}y_n),
$$
where $\epsilon$ is a primitive $a_i$th root of unity.
The identification is given by $y_j^{a_i}=x_j^{a_i}/x_i^{a_j}$. 
(\ref{wps.defn}.1) are called the {\it orbifold charts} on
 $\p(a_0,\dots,a_n)$.

For any $i$, 
$\p(a_0,\dots,a_n)$ has  an index $a_i$ quotient singularity at $P_i$.
For any $i<j$, if
$\gcd(a_0,\dots, \hat{a_i},\dots,\hat{a_j},\dots,a_n)>1$,
then $\p(a_0,\dots,a_n)$ has  a  quotient singularity 
along $(x_i=x_j=0)$.
These give all the codimension 2 singular subsets of
$\p(a_0,\dots,a_n)$.

For every $m\in \z$  there is a  rank 1 sheaf
$\o_{\p}(m)$ which is locally free only if $a_i|m$ for every $i$.
A basis of the space of  sections of $\o_{\p}(m)$
is given by all monomials in $x_0,\dots ,x_n$ with weighted degree  $m$.
Thus $\o_{\p}(m)$ may have no sections for some $m>0$.
\end{defn}

\begin{say}[Anticanonically embedded quasi smooth Fano hypersurfaces]{\ }
\label{ant.emb.say}

Let $X\in |\o_{\p}(m)|$ be a hypersurface of degree $m$. 
The adjunction formula
$$
K_X\cong \o_{\p}(K_{\p}+X)|_X \cong \o_{\p}(m-(a_0+\cdots +a_n))|_X
$$
holds if  $X$ does not contain any of the codimension 2 singular subsets.
If this condition is satisfied then $X$ is a  Fano variety
iff $m<a_0+\cdots +a_n$.  Frequently
the most interesting
cases are when $m$ is as large as possible. Thus we consider the case
$X_d\in |\o_{\p}(d)|$ for $d=a_0+\cdots +a_n-1$.
Such an $X$ is called {\it anticanonically embedded}.

In most  cases, all  hypersurfaces  of a given degree $d$
are singular  and   pass through
some of the vertices $P_i$. In these cases  the best one can hope 
is that a general hypersurface $X_d$  is smooth in the
orbifold sense, called {\it quasi smooth}. At the vertex $P_i$ this means that
the preimage of $X_d$ in the orbifold chart 
$\c^n(y_0,\dots, \hat{y_i},\dots,y_n)$
is smooth. In terms of the monomials of degree $d$  this is equivalent 
to saying that
$$
\mbox{For every $i$ there is a $j$ and 
 a monomial $x_i^{m_i}x_j$ of degree $d$.}
\eqno{(\ref{ant.emb.say}.1)}
$$
$j=i$ is allowed, corresponding to the case when
the general $X_d$ does not pass through $P_i$. 
The condition that $X_d$ does not contain any of the singular 
codimension 2 subsets
is equivalent to
$$
\begin{array}{l}
\mbox{If $\gcd(a_0,\dots, \hat{a_i},\dots,\hat{a_j},\dots,a_n)>1$
 then there is a}\\  
\mbox{monomial of degree $d$ not involving $x_i, x_j$.}
\end{array}
\eqno{(\ref{ant.emb.say}.2)}
$$
For $n\geq 3$, these are the two most important special cases of the general
quasi smoothness condition:
$$
\begin{array}{l}
\mbox{For every $I\subset \{0,\dots,n\}$  there  is an injection}\\
\mbox{ $e: I\hookrightarrow \{0,\dots,n\}$ and monomials}\\
\mbox{ $x_{e(i)}\prod_{j\in I}x_j^{m_{ij}}$ 
of degree $d$ for every $i\in I$.}
\end{array}
\eqno{(\ref{ant.emb.say}.3)}
$$
\end{say}

\begin{rem} The quasi-smoothness condition in \cite[I.5.1]{fletch}
says that 
$$
\begin{array}{l}
\mbox{For every $I\subset \{0,\dots,n\}$ either (\ref{ant.emb.say}.3) holds,}\\
\mbox{or there is a monomial $\prod_{j\in I}x_j^{b_{j}}$ 
of degree $d$.}\\
\end{array}
\eqno{(\ref{ant.emb.say}.3')}
$$
 The two versions are, however, equivalent.
We prove this by induction on $|I|$.
Indeed, assume that there is a monomial
$\prod_{j\in I}x_j^{b_{j}}$ 
of degree $d$ and let $I'\subset I$ be all the indices which are involved in
at least one such monomial. By induction (\ref{ant.emb.say}.3)
holds for $I\setminus I'$, giving monomials
$x_{e(i)}\prod_{j\in I\setminus I'}x_j^{b_{ij}}$ for $i\in  I\setminus I'$. 
By assumption these $e(i)$ are not in $I$, so we can choose
$I'\cup e(I\setminus I')$ as the image of $e:I\to \{0,\dots,n\}$.
(A suitable reordering of the values of $e$ may be necessary.)
\end{rem}

The computer searches done in connection with \cite{fletch} and 
\cite{dk} looked at values of $a_i$ in a certain range
 to find the $a_i$ satisfying the constraints
(\ref{ant.emb.say}.1-3).  This approach starts with the $a_i$ and views
(\ref{ant.emb.say}.1-3) as linear equations in the unknown 
nonnegative integers $m_i, m_{ij}$.

In the cases studied in \cite{fletch} and 
\cite{dk} these searches seemed exhaustive. Aside from one series of
examples, the computers produced solutions for low values of the $a_i$
and then did not find any more as the range of the allowable values 
was extended. This of course does not ever lead to a proof that the
lists were complete.

A similar search for log Fano hypersurfaces in weighted projective
4--spaces seems much harder. With some reasonably large bounds,
say $a_i\leq 100$, the programs run very slowly and they produce  examples
where the $a_i$ are near 100. It is quite unlikely that
 any systematic search of this kind could have discovered the   example
with the largest $a_0$:

$(  407,   547,  5311, 12528, 18792)$ with monomials

$$x_4^{2},\ x_3^{3},\ x_1^{59}x_2,\ x_0x_2^{7},\ x_0^{91}x_1, $$

\noindent or the beautiful pair of sporadic examples with largest $a_4$:

$(  223,  9101, 46837, 112320, 168480)$ with monomials

$$x_4^{2},\ x_3^{3},\ x_1x_2^{7},\ x_0x_1^{37},\ x_0^{1301}x_2,\qtq{and} $$

$(  253,  7807, 48101, 112320, 168480)$ with monomials

$$x_4^{2},\ x_3^{3},\ x_1^{37}x_2,\ x_0x_2^{7},\ x_0^{1301}x_1. $$

\noindent The biggest values of $a_0$ are of some interest in connection
with the conjectures of \cite[1.3]{shok}.

Next we describe the computer programs that led to the
list of  anticanonically embedded quasi smooth 
Fano hypersurfaces in weighted projective 4-spaces.
The  programs, written in C,  are available at
\begin{verbatim}www.math.princeton.edu/~jmjohnso
\end{verbatim}

\begin{say}[Preliminary steps]{\ }\label{descr.of.progr}

In order to find all solutions, we change the point of view.  We
consider 
(\ref{ant.emb.say}.1) to be the main constraint with 
  coefficients $m_i$ and 
unknowns $a_i$. The corresponding equations  can then be written
as a linear system 
$$
(M+J+U) (a_0\ a_1\ a_2\ a_3\ a_4)^t
= (\ -1\ -1\ -1\ -1\ -1)^t
\eqno{(\ref{descr.of.progr}.1)}
$$
where $M=\diag(m_0,m_1,m_2,m_3,m_4)$ is a diagonal matrix,
$J$ is a matrix with all entries $-1$ and $U$ is
a matrix where each row has 4 entries $=0$ and one entry $=1$.
The main advantage is that
some of the $m_i$ can be bounded a priori.
Assume for simplicity that $a_0\le a_1\le \cdots \le a_4$.

Consider for instance $m_4$. The relevant equation is
$$
m_4a_4+a_{e(4)}=a_0+a_1+a_2+a_3+a_4-1.
$$
Since $a_4$ is the biggest, we get right away that $1\leq m_3\leq 3$.
 Arguing inductively
with some case analysis we obtain that
$$
3\leq m_2\leq 16,\ 2\leq m_3\leq 6,\ 1\leq m_4\leq 3
\eqno{(\ref{descr.of.progr}.2)}
$$

Thus we have only finitely many possibilities for
the matrix $U$ and the numbers $m_2,m_3,m_4$.
Fixing these values, we obtain a linear system 
$$
(M+J+U) (a_0\ a_1\ a_2\ a_3\ a_4)^t
= (\ -1\ -1\ -1\ -1\ -1)^t,
$$
where the only variable coefficients are $m_0,m_1$
in  the upper left corner
of $M$. Solving these formally we obtain that
$$
a_0=\frac{\alpha_0m_1+\beta_0}{\gamma_2m_0m_1+\gamma_0m_0+\gamma_1m_1+\delta},
a_1=\frac{\alpha_1m_0+\beta_1}{\gamma_2m_0m_1+\gamma_0m_0+\gamma_1m_1+\delta}.
$$
where the $\alpha_i,\beta_i,\gamma_i,\delta$
 depend only on $U$ and  $m_2,m_3,m_4$.
\medskip

We distinguish 3 cases. The first one is the main source of
examples. Cases 2 and 3 are anomalies from the point of view of our method.
In both cases we ended up experimentally finding strong restrictions
on the $a_i$. Even with hindsight we do not know how to prove these a priori.

\begin{case}[$\gamma_2\neq 0$] In this case the absolute value of
$$
\frac{\alpha_0m_1+\beta_0}{\gamma_2m_0m_1+\gamma_0m_0+\gamma_1m_1+\delta}
$$
goes to zero as $m_0,m_1$ go to infinity. It is not hard to write
 down the precise condition and a computer check shows that
$$
\frac{\alpha_0m_1+\beta_0}{\gamma_2m_0m_1+\gamma_0m_0+\gamma_1m_1+\delta}\geq 1
\quad \Rightarrow \quad   \min\{m_0,m_1\}\leq 83.
$$
\end{case}

\begin{case}[$\gamma_2=0,\gamma_0\gamma_1\neq 0$ ] 
 It turns out that if this holds then $\gamma_0\gamma_1>0$
and $a_0,a_1$ are bounded by 8 for $\min\{m_0,m_1\}\geq 36$. 
Moreover, the 3 linear forms 
$$
\alpha_0m_1+\beta_0,\ \alpha_1m_0+\beta_1,\ \gamma_0m_0+\gamma_1m_1+\delta
$$
are dependent. This implies that
$$
\alpha_1\gamma_1a_0+\alpha_0\gamma_0a_1=\alpha_0\alpha_1.
$$
A computer search shows that this is possible only if $a_0=a_1=1$.
\end{case}

\begin{case}[$\gamma_2=0,\gamma_0\gamma_1=0$] 
It turns out that  one of 
$\alpha_0m_1+\beta_0, \alpha_1m_0+\beta_1$
 equals $\gamma_0m_0+\gamma_1m_1+\delta$. Thus $a_0=1$ or $a_1=1$.
Moreover, we also see by explicit computation that
one of the following holds:
$$
a_2=a_3=a_4, \quad  a_2=a_3=a_4/2 \qtq{or} a_2=a_3/2=a_4/3.
$$

\end{case}
\end{say}

\begin{say}[Main Computer Search]{\ }\label{main.series}

Here we discuss the main case when, in addition to the inequalities
(\ref{descr.of.progr}.2) we also assume that $3\leq m_1\leq 83$.
In this case the system (\ref{descr.of.progr}.1)
reduces to a single unknown $m_0$. This is very similar to the
4--variable case discussed in \cite{jk1}.

We solve formally for $a_0$ to get
$$
a_0=\frac{\gamma_0}{m_0\alpha+\beta}
$$
where $\alpha,\beta,\gamma_0$ depend only on $U$ and  $m_1,m_2,m_3,m_4$.
If
$\alpha\neq 0$ then we get a bound on $m_0$ too, and we are down to finitely
many possibilities all together.
There are 403455  cases of this.
The resulting solutions need considerable cleaning up.
Many of them  occur multiply and we also have to check
the other conditions (\ref{ant.emb.say}.2-3).  Discarding repetitions, we get
15757 cases, out of which 4594 are quasi smooth.

 If $\alpha=0$ then
we get a series solution where the $a_i$ are linear functions of a variable
 $m_0$. There are 550122  cases of this.
Here the main difficulty is that the program does not produce the series
in a neat form. Usually one series is put together out of
many pieces according to some congruence condition.
\end{say}

\begin{say}[Additional  Cases]{\ }\label{a0=1.series}

Assume first that we are in Case 2 of (\ref{descr.of.progr}).
Since $a_0=a_1=1$, the numbers $a_1,a_2,a_3,a_4$ and 
  $d=a_1+a_2+a_3+a_4$ satisfy the numerical conditions
(\ref{ant.emb.say}.3). 
 This leads to a
lower dimensional problem which is easy to solve.

Case 2 of (\ref{descr.of.progr}) is even easier. We get
solutions of the form
$$
(1,a,b,b,b),\quad (1,a,b,b,2b) \qtq{or} (1,a,b,2b,3b).
$$
Applying (\ref{ant.emb.say}.3) to $I=\{2,3,4\}$
gives that $b|a$. Thus $b$ divides all but one of the weights,
so $b=1$. This implies that $a\leq 6$.
At any case, all these appear also under Case 2 of (\ref{descr.of.progr}).
\end{say}

At the end we obtain our first main result:

\begin{thm}\label{list.thm}
 The following is a complete list of 
anticanonically embedded quasi smooth 
Fano hypersurfaces in weighted projective 4-spaces:
\begin{enumerate}
\item 48 infinite series of the form
$$
X_{2k(b_1+b_2+b_3)}\subset
\p(2,kb_1,kb_2,kb_3,k(b_1+b_2+b_3)-1)\qtq{for $k=1,3,5,\dots$}
$$
The occurring 3-tuples $b_1,b_2,b_3$ are described in (\ref{3-tuple.rem}).
\item 4442 sporadic examples whose list is available at
\begin{verbatim}www.math.princeton.edu/~jmjohnso
\end{verbatim}
\end{enumerate}
\end{thm}

\begin{say}[An  error check]{\ }

We wrote a program that looked at all 5-tuples satisfying
$$
a_0\leq 100,\ a_1\leq 200,\ a_2\leq 200,\ a_3\leq 400,\ a_4\leq 600.
$$
The program  ran for 4 days and produced
3610 quasi smooth examples.
These were in complete agreement with the
correspondingly truncated  list of 4442 sporadic examples.
\end{say}

\begin{rem}\label{3-tuple.rem} It turns out that a 
3-tuple $b_1,b_2,b_3$ appears in (\ref{list.thm}.1)
iff $|-2K|$ of $\p(b_1,b_2,b_3)$ has a quasi smooth member. 
The list of these is implicit in Reid's list of 95 families of singular K3
surfaces in weighted projective 3-spaces. 
In \cite[II.3.3]{fletch} they correspond to those quadruplets
$(b_1,b_2,b_3,b_4)$ for which $b_4=b_1+b_2+b_3$.
Our 48 3-tuples occured
explicitly in \cite{yonemura, tomari} in connection with
the study of simple K3 singularities of multiplicity 2.

One direction of this observation is easy to establish in all dimensions.

\begin{lem} Assume that $|-2K|$ of $\p(b_1,\dots,b_n)$ has a 
quasi smooth member.  Then the general 
anticanonically embedded 
Fano hypersurface in
$$
\p(2,kb_1,\dots,kb_n,k(b_1+\dots+b_n)-1)
$$
is quasi smooth for $k=1,3,5,\dots$
\end{lem}

We conjecture that  conversely, every infinite series is
of this form. It is interesting that every 
quasi smooth
hypersurface
in  $\p(2,kb_1,\dots,kb_n,k(b_1+\dots+b_n)-1)$
has a codimension 2 singular set. Thus the above conjecture would imply
that for every $n\geq 4$ there are only finitely many
anticanonically embedded quasi smooth 
Fano hypersurfaces with isolated singularities 
in weighted projective $n$-spaces.
\end{rem}

It is not hard to check  which of the above 
Fano threefolds have terminal singularities. 
The families in (\ref{list.thm}.1) always have nonisolated singularities,
and for the remaining cases the conditions of \cite[II.4.1]{fletch}
work. 
As a consequence, 
 we obtain the following:

\begin{cor}\label{95.cor} 
 The Reid--Fletcher list (cf.\ \cite[II.6.6]{fletch})  of
95  families of anticanonically embedded quasi smooth 
terminal Fano threefolds
in  weighted projective 4-spaces is complete. \qed
\end{cor}

Next we study the existence of K\"ahler--Einstein metrics and
 the nonexistence of tigers on our Fano hypersurfaces.
After some definitions we recall the criterion established in
\cite{jk1}. In the case of K\"ahler--Einstein metrics
this in turn relies on earlier work of \cite{nadel-KE, dk}.

\begin{defn}\label{klt.etc.defn}
 Let $X$ be a normal variety and $D$ a $\q$-divisor on $X$.
Assume for simplicity that $K_X$ and $D$ are both $\q$-Cartier.
Let $g:Y\to X$ be any proper birational morphism, $Y$ smooth.
 Then there is a unique
$\q$-divisor $D_Y=\sum e_iE_i$ on $Y$ such that 
$$
K_Y+D_Y\equiv g^*(K_X+D)\qtq{and}  g_*D_Y=D.
$$
We say that $(X,D)$ is
{\it klt} (resp.\ {\it log canonical}) if
$e_i>-1$ (resp.\ $e_i\geq -1$) for every $g$
and for every $i$.
See, for instance,  \cite[2.3]{kmbook} for a detailed  introduction.
\end{defn}

\begin{defn} \cite{keel-mc}\label{tiger.defn}
 Let $X$ be a normal variety. A {\it tiger}
on $X$  is an effective  $\q$-divisor $D$ such that $D\equiv -K_X$ 
and $(X,D)$ is not klt.
As illustrated in \cite{keel-mc}, the tigers  carry
 important  information about birational transformations of 
log del Pezzo surfaces. They are expected to play a similar role
in higher dimensions.
\end{defn}

\begin{prop}\cite{jk1} \label{main.estimate.cor}
  Let $X_d\subset \p(a_0,\dots,a_n)$ be a quasi smooth
hypersurface of degree $d=a_0+\cdots +a_n-1$. 
\begin{enumerate}
\item $X$ does not have a tiger
if $d\leq a_0a_1$.
\item  $X$ admits a
K\"ahler--Einstein metric if
 $d< \tfrac{n}{n-1} a_0a_1$.\qed
\end{enumerate}
\end{prop}

\begin{cor} \label{KE.cor}
Of the sporadic series of quasi smooth Fano hypersurfaces
given in (\ref{list.thm}.2),
there are 1605 types where none of the members have a tiger
and 1936 types where every  member admits a K\"ahler--Einstein metric.
This information is contained in the list given in (\ref{list.thm}.2).
\end{cor}

Finally we study the case of 
 Calabi--Yau hypersurfaces 
and hypersurfaces  of general type in 
weighted
projective spaces. For these cases  there are finiteness results in
all dimensions. The key part is the case of 
 Calabi--Yau hypersurfaces.

\begin{thm}\label{CY.thm} For any $n$ there are only finitely many types
of quasi smooth hypersurfaces with trivial canonical class in
weighted projective spaces $\p(a_0,\dots,a_n)$.
\end{thm}

Proof. As in the Fano case, first we look at those hypersurfaces
which are quasi smooth at the vertices of $\p(a_0,\dots,a_n)$.
This condition is equivalent to a linear system of
equations
$$
(M+J+U) (a_0, \dots, a_n)^t
= (0, \dots, 0)^t
\eqno{(\ref{CY.thm}.1)}
$$
where $M=\diag(m_0,\dots,m_n)$ is a diagonal matrix,
$J$ is a matrix with all entries $-1$ and $U$ is
a matrix where each row has $n$ entries $=0$ and one entry $=1$.
In the geometric setting  the $m_i$ and the $a_i$ are positive integers,
but it will be convenient to allow the $a_i$ to be positive real numbers.
By the homogenity of the system we may assume that $\sum a_i=1$.

Assume now that we have an infinite sequence of solutions
$(a_0(t), \dots, a_n(t))$ where a priori $M(t), J(t), U(t)$
also vary with $t$. By passing to a subsequence
we may assume that $J(t),U(t)$ are constant and each $a_i(t)$ converges to 
a value $A_i$.  Thus we can write
$a_i(t)=A_i+c_i(t)$ where $\lim_{t\to\infty} c_i(t)=0$,
$\sum_i A_i=1$ and $\sum_i c_i(t)=0$.
By passing to a subsequence and rearranging, we can also assume that
$I:=\{i: c_i(t)<0\}$ is independent of $t$ and that
$A_0/(-c_0(t))$ is the smallest positive number among
$\{A_i/(-c_i(t)): i\in I\}$. The quasi smoothness condition at the vertex $P_0$
translates into $m_0(t)a_0(t)+a_j(t)=1$.
$\lim_{t\to\infty}a_0(t)=A_0>0$ since $c_0(t)<0$, hence 
$m_0(t)$ is bounded from above.
Thus we may assume that $m_0(t)=m_0$ is constant and 
$\lim_{t\to\infty}m_0c_0(t)+c_j(t)=0$.
$m_0a_0(t)+a_j(t)=1$
 is equivalent to
$$
[m_0A_0+A_j]+[m_0c_0(t)+c_j(t)]=1.
\eqno{(\ref{CY.thm}.2)}
$$
By the above considerations, 
(\ref{CY.thm}.2) splits into 2 equations
$$
m_0A_0+A_j=1 \qtq{and} m_0c_0(t)+c_j(t)=0.\eqno{(\ref{CY.thm}.3)}
$$
Using $\sum_i c_i(t)=0$ and the second equation in (\ref{CY.thm}.3)
we obtain that
$$
\sum_{i\in I}c_i(t)=-\sum_{i\not\in I}c_i(t)\leq -c_j(t)=m_0c_0(t).
\eqno{(\ref{CY.thm}.4)}
$$
Multiplying by $A_0/c_0(t)$ and using the special choice of $A_0/c_0(t)$
we get that
$$
m_0A_0\leq \sum_{i\in I}c_i(t)\frac{A_0}{c_0(t)}\leq \sum_{i\in I} A_i.
\eqno{(\ref{CY.thm}.5)}
$$
Combining with the first equation of (\ref{CY.thm}.3) we get that
$$
1=m_0A_0+A_j\leq A_j+\sum_{i\in I} A_i\leq \sum_{i=0}^n A_i=1.
\eqno{(\ref{CY.thm}.6)}
$$
This implies that all inequalities in (\ref{CY.thm}.4--6)
are equalities. Hence $A_k,c_k(t)$ are zero for $k\not\in I\cup\{j\}$.
By assumption the $a_k(t)$ are positive, so $I\cup\{j\}=\{0,\dots, n\}$. 
Moreover, the ratios $A_i/c_i(t)$ are all the same for $i\in I$. 

These imply that, up to rearranging the indices, the 
 $a_i(t)$ are of the form
$$
(A_0(1-c(t)),\dots, A_{n-1}(1-c(t)), A_n+c(t)\textstyle{\sum_{i=0}^{n-1}}A_i).
$$
Consider next the equation
$$
m_n(A_n+c(t)\textstyle{\sum_{i=0}^{n-1}}A_i)+A_j(1-c(t))=1,
$$
where for notational simplicity we allow $j=-1$ with $A_{-1}=0$.
For large $t$ this implies that 
$\textstyle{\sum_{i=0}^{n-1}}A_i=A_j$, which is not possible for $n\geq 2$. 
Thus  $A_n=0$ and the solutions become
$$
(A_0(1-c(t)),\dots, A_{n-1}(1-c(t)), c(t))\qtq{where} 
\textstyle{\sum_{i=0}^{n-1}}A_i=1.
\eqno{(\ref{CY.thm}.7)}
$$
To get quasi smoothness, we need to understand all monomials of
degree $\sum a_i$, which amounts to finding all integer solutions of
$\sum b_ia_i=1$.  In our case, for large $t$ there are no solutions
with $b_n=0$ which means that  every hypersurface
of degree $\sum a_i$ contains the hyperplane $(x_n=0)$, hence
they are all reducible. Thus the solutions
(\ref{CY.thm}.7) do not correspond to
quasi smooth hypersurfaces.\qed

\begin{rem} The solutions (\ref{CY.thm}.7) do correspond to
interesting series of singularities. Namely, 
for every integer solution of $\textstyle{\sum_{i=0}^{n-1}}\tfrac1{m_i}=1$
they give an infinite series of singularities
$$
(x_0^{m_0}+\cdots +x_{n-1}^{m_{n-1}}+x_n^k)x_n=0\subset \a^{n+1}
\qtq{for $k=1,2,\dots$.}
$$
These singularities are weighted homogeneous and
semi log canonical  (see \cite[16.2.1]{ketal} for the definition)
but not isolated.
By adding a general higher degree term, we get isolated log canonical
singularities.
\end{rem}

\begin{cor}\label{gen.t.cor} For any $n$ and $k >0$ 
there are only finitely many families
of quasi smooth hypersurfaces  $X\subset \p(a_0,\dots,a_n)$
such that $\omega_X\cong \o_X(k)$. 
\end{cor}

Proof. Assume that
$$
X=(F(x_0,\dots,x_n)=0)\subset \p(a_0,\dots,a_n)
$$
 is quasi smooth of degree $d$ and
$\omega_X\cong \o_X(k)$. Then
$$
X^*:=(F(x_0,\dots,x_n)+x_{n+1}^d+\dots+x_{n+k}^d=0)
\subset \p(a_0,\dots,a_n,\overbrace{1,\dots,1}^{k-times})
$$
is also quasi smooth of degree $d$ and
$\omega_X\cong \o_X$. Thus we are done by (\ref{CY.thm}).\qed

\begin{rem} The finiteness result (\ref{gen.t.cor})
is in accordance with the conjectures \cite[18.16]{ketal}.
On the other hand, (\ref{CY.thm}) seems to be a more special
finiteness assertion.
\end{rem}

\begin{ack}  We   thank  J.\ McKernan
for  helpful
comments and references.
Partial financial support was provided by  the NSF under grant number 
DMS-0096268. 
\end{ack}

\vskip1cm

\noindent Princeton University, Princeton NJ 08544-1000

\begin{verbatim}jmjohnso@math.princeton.edu\end{verbatim}

\begin{verbatim}kollar@math.princeton.edu\end{verbatim}
\end{document}